\UseRawInputEncoding
\documentclass[12pt]{article}

\usepackage{hyperref}
\usepackage{booktabs}

\usepackage[numbered,framed]{matlab-prettifier}

\usepackage{enumitem}

\usepackage{hyperref} \usepackage{breakurl}

\usepackage{listings}
\usepackage{adjustbox}
\usepackage[margin=1in,footskip=0.25in]{geometry}
\usepackage{relsize}
\usepackage{amsthm}%
\usepackage{mathtools}
\usepackage{amsfonts}%
\usepackage{pdfpages}
\usepackage{amsmath}
\usepackage{color}
\usepackage{flexisym}
\definecolor{mygreen}{RGB}{28,172,0} 
\definecolor{mylilas}{RGB}{170,55,241}

\usepackage{cite}
\newtheorem{theorem}{Theorem}
\newtheorem{remark}[theorem]{Remark}

\newtheorem{conjecture}[theorem]{Conjecture}

\begin{document}

\title{On the stability domain of a class of linear systems of fractional order}

\author{Marius-F. Danca{\footnote{Corresponding author}}\\
STAR-UBB Institute,\\
Babes-Bolyai University,\\
400084, Cluj-Napoca, Romania\\
and\\
Romanian Institute of Science and Technology, \\
400487, Cluj-Napoca, Romania\\
Email: danca@rist.ro\\
}

\maketitle
\begin{abstract}
In this paper, the shape of the stability domain $S^q$ for a class of difference systems defined
by the Caputo forward difference operator $D^q$ of order $q\in (0, 1)$ is numerically analyzed. It is
shown numerically that due to of power of the negative base in the expression of the stability domain,
in addition to the known cardioid-like shapes, $S^q$ could present supplementary regions where the
stability is not verified. The Mandelbrot map of fractional order is considered as an illustrative
example. In addition, it is conjectured that for $q < 0.5$, the shape of $S^q$ does not cover the main body
of the underlying Mandelbrot set of fractional order as in the case of integer order.
\end{abstract}


\vspace{5mm}

%
%

\section{Introduction}

In the last decade, fractional calculus has been considered as an important tool in many scientific and engineering fields.
The basic theory of fractional calculus modeling and control
systems can be found in, e.g.,  \cite{xy1,xy2,xy3,xy4,xy5,xy6,xy8}. Studies of their applications are presented in \cite{yy1,yy2,yy3}. Systems with fractional variable orders are analyzed in \cite{zz1,zz2}.
In terms of the stability analysis of fractional differential equations, one of the important properties that is analysed in
order to study the behavior of the considered systems is presented in, e.g.,  \cite{tt1,tt2,tt3,bibsapte},  {while applications and stability studies of discrete fractional difference equations can be found in \cite{x1,x2,x3}}.

For commensurate fractional-order systems, several powerful criteria are established. The most well-known Matignon's stability theorem \cite{tt5} (see also \cite{xy7}) determines the system stability by searching for the location of the eigenvalues in the complex plane, which represents a starting point for most research into the stability of fractional-order systems.

In \cite{bibunu}, a stability criterion for a fractional difference linear system is presented:
\begin{equation}\label{unu}
\Delta^q y(n+1-q)=Ay(n),~~n=0,1,...,
\end{equation}
where $A\in \mathbb{R}^{d\times d}$ and $\Delta ^q$ is the Caputo forward difference operator (see, e.g.,  \cite{bibpatru,bibcinci}).

As shown in \cite{bibunu}, the stability criterion for \eqref{unu} is proved to be fully explicit, also involving the decay rate of the solutions.

The asymptotic stability property of \eqref{unu} for both scalar and vector cases is stated by the following result:

\begin{theorem}[\cite{bibdoi}] \label{t1}
Let $q\in(0,1)$ and $A\in \mathbb{R}^{d\times d}$. Then, \eqref{unu} is asymptotically stable if and only if the isolated zeros, off the non-negative real axis, of $\det (I-z^{-1}(1-z^{-1})^{-q}A)$ lie inside the unit circle.
\end{theorem}

In \cite {bibunu}, the aim is to formulate an alternative stability criterion for the integer difference system:
\begin{equation}\label{doi}
\Delta y(n)=Ay(n),~~n=0,1,...,
\end{equation}
where $\Delta y(n)=y(n+1)-y(n)$ is the standard operator here for $q=1$.

As is well known, the system \eqref{doi} is asymptotically stable if and only if all the eigenvalues of $I + A$, where $I$ is the identity matrix, are located inside the unit circle. If one considers the stability set for Equation \eqref{doi} in the polar form

\begin{equation*}
S=\bigg\{z\in C: |z|<-\cos(\arg z)~~ \text{and}~~~ |\arg z|>\frac{\pi}{2}\bigg\},
\end{equation*}
the stability result can be reformulated in the following form:
\begin{theorem}[\cite{bibunu}]
The linear difference system \eqref{doi} is asymptotically
stable if and only if $\lambda\in S$ for all the eigenvalues $\lambda$ of $A$. In this case, the solutions of \eqref{doi} decay towards zero exponentially as $n\rightarrow \infty$.
\end{theorem}

In order to introduce an alternative stability criterion for \eqref{unu}, as a direct extension of Theorem \eqref{t1}, consider the following set:

\begin{equation}\label{main}
S^q=\Bigg\{z\in C: |z|<\left(2\cos\frac{|\arg z|-\pi}{2-q}\right)^q~~ \text{and}~~~ |\arg z|>\frac{q\pi}{2}\Bigg\},
\end{equation}
Note that the second inequality in \eqref{main} represents the Matignon criterion.

The main result in \cite{bibunu} is:

\begin{theorem}[\cite{bibunu}]\label{x}
Let $q\in(0, 1)$ and $A\in \mathbb{R}^{d\times d}$. If $\lambda\in S^q$ for all the eigenvalues $\lambda$ of $A$, then system \eqref{unu} is asymptotically stable. In this case, the
solutions of \eqref{unu} decay towards zero in such a way that
\[
\|y(n)\| = O(n^{-q})~~ as~~ n\rightarrow \infty,
\]
for any solution $y$ of \eqref{unu}. Furthermore, if $\lambda\in C\setminus cl(S^q)$ for an eigenvalue $\lambda$ of $A$, then \eqref{unu} is
not stable.
\end{theorem}

Because of the explicit and convenient form of $S^q$, drawn in the complex plane for eigenvalues, Theorem \ref{x} became widely used in practical applications.

 In this paper, it is shown that for some empirically found values of $q\in(0,1)$, the form of $S^q$ might present other unexpected shapes. In addition, a conjecture related to the shape of $S^q$ in the case of the Mandelbrot set of fractional order is introduced.

\section{On the Shape of the Stability Area \boldmath{$S^q$}}\label{sec2}
The parametric equations of the frontier of the stability region $S^q$ for a fixed point, $z=x+\i y$, $x,y\in \mathbb{R}$, can be easily drawn in the following form \cite{bibunu}:
\begin{equation}\label{para}
\begin{array}{l}
x=-2^q\cos^q\theta \cos((2-q)\theta)\\
y=-2^q\cos^q\theta \sin((2-q)\theta),~~ |\theta|\leq \pi/2.
\end{array}
\end{equation}

On the other hand, to analyze the shape of $S^q$ in more detail, consider the form \eqref{main}.
As is known, the argument of a complex number $z=x+iy$, $\arg z$, used in \eqref{main}, is the angle between the positive real axis and the line joining the origin and the image in the complex plane of $z$.
Usually, the value of the argument is numerically computed as
\begin{equation*}\label{cinci}
\text{arg} z=\text{atan2}(y,x),
\end{equation*}
where the two-argument function $\text{atan2}$ is an available implemented function in the math libraries of many programming languages.

In contrast to the inverse tangent function, $\tan ^{-1}$ ($\arctan$, or atan), the function $\text{atan2}$ computes the principal value of the arctangent of $y/x$, determining the quadrant of the returned value by using the signs of both arguments. Note that a domain error may occur if both arguments are zero. (Regarding the 1999 ISO C, 1978 ANSI Fortran, or 1982 ISO Pascal Standard standard descriptions and other issues related to $\text{atan2}$ for the case of $x=0$ and $y=0$, they are referred to on p. 70 \cite{bibtrei}).
However, as will be shown below, the problem of arg in \eqref{main} can be quite complicated, even if both arguments are not zero, as shown next with simple reasoning.

While $\arctan$ gives an angle between $-\pi/2$ and $\pi/2$, the function $\text{atan2}$ always provides a result within $(-\pi,\pi)$.

\begin{figure}[h]
\includegraphics[scale=1.2]{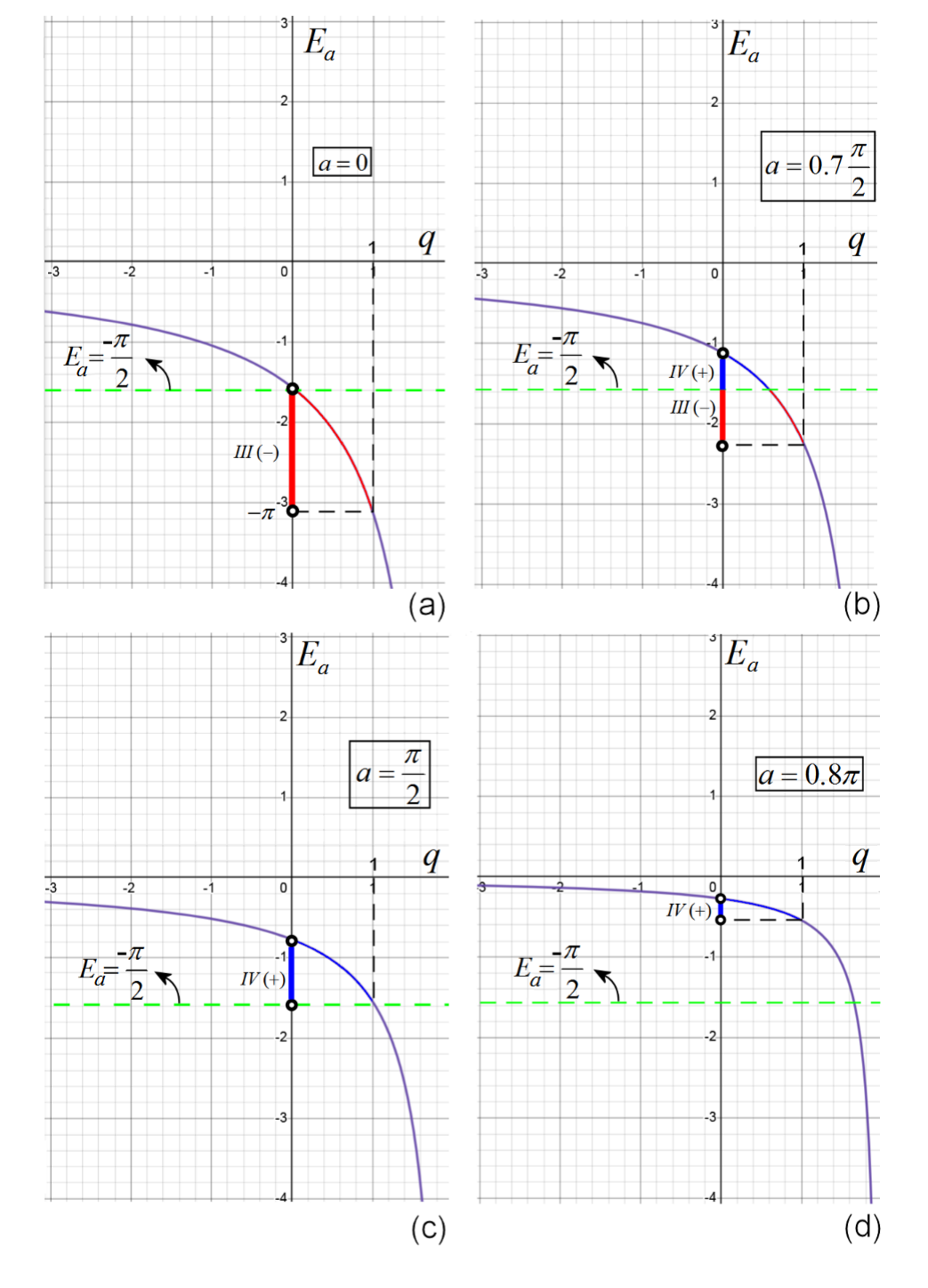}
\caption{Graph of the function $E_a$, for four representative cases: (\textbf{a}) the limit case $a=0$; \mbox{(\textbf{b}) $a=0.7\pi/2$;} (\textbf{c}) $a=\pi/2$; and (\textbf{d}) $a=0.8\pi$.}
\label{fig1}
\end{figure}

Denote, in \eqref{main},
\begin{equation}\label{a}
a:=|\arg z|,
\end{equation}
and consider the argument of the cos function in $S^q$ as the function $E_a:\mathbb{R}\setminus \{2\}\rightarrow \mathbb{R}$, defined as
\begin{equation*}\label{e}
E_a(q)=\frac{a-\pi}{2-q},
\end{equation*}
with $a$ being a parameter (Figure \ref{fig1}). Being determined with atan2, $a$ given by \eqref{a} satisfies the following relations: $a\in[0,\pi)$, and $a-\pi\in[-\pi,0)$.
Therefore, because for $q\in(0,1)$, $2-q>0$, one has

\[
E_a(q)<0.
\]

In order to determine in which quadrant the values of $E_a(q)$ are situated, let us analyze and sketch the graph of the restriction of the function $E_a$ to $(0,1)$. The derivative of $E_a$ is $E'_a(q)=(a-\pi)/(2-q)^2<0$. Therefore, from the monotonicity, one deduces that $E_a$ decreases from $(a-\pi)/2$, for $q\downarrow 0$, but to $a-\pi$, for $q\uparrow1$. For example, for the limit case of $a=0$, $E_0(q)$ takes values within the range $(-\pi,-\pi/2)$, i.e., $E_0(q)$ are situated in quadrant III (Figure  \ref{fig1}a). For the case of $a=0.7\pi/2$, the values of $E_{\scriptscriptstyle{0.7}\frac{\pi}{2}}(q)$ can be situated in both quadrants III and IV (Figure  \ref{fig1}b). For $a=\pi/2$, $E_{\frac{\pi}{2}}$ moves within in the quadrant IV (Figure  \ref{fig1}c). The case of $a=0.8\pi$ is presented in Figure  \ref{fig1}d.
Therefore, for $q\in(0,1)$ and $a\in[0,\pi)$, the argument $E_a$ can take values situated in quadrant III and/or IV, where $\cos(E)$ could be either negative or positive.

The relation between $q$ and $a$ on the boundary between the two quadrants III and IV is
\begin{equation*}\label{x1}
q=2a/\pi.
\end{equation*}

\begin{figure}
\centering
\includegraphics[scale=1.55]{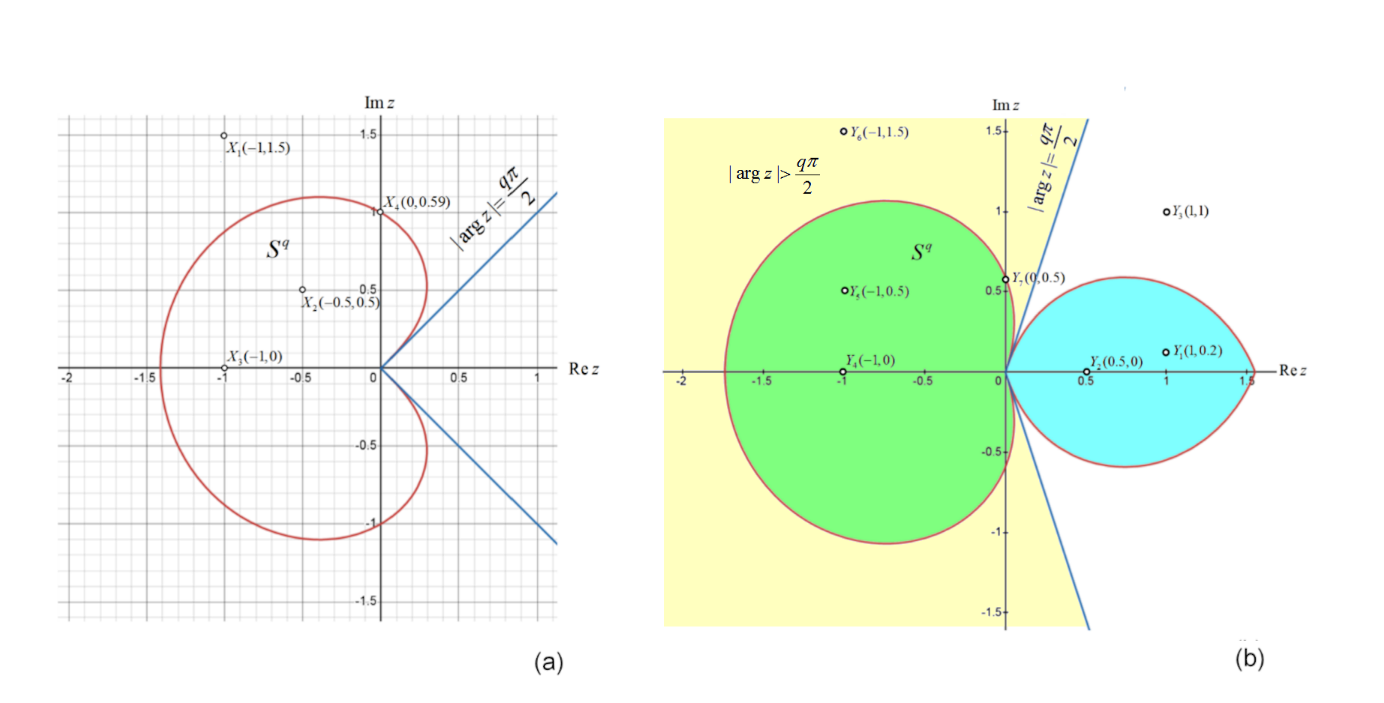}
\caption{Graph of $\Gamma$ (red plot) for two representative cases: (\textbf{a}) $q=0.5$; (\textbf{b}) $q=0.8$. Green plot represents the stability domain $S^q$, surrounded by $\Gamma$, delimitated by the inequality $|\arg z|>q\pi/2$, while light blue and white domains do not belong to the stability domain.}
\label{fig2}
\end{figure}

To obtain the shape of $S^q$ via relation \eqref{main}, consider a variable fixed point with the underlying variable eigenvalue, $z=x+iy$, to be replaced in \eqref{main}. Denote hereafter by $\Gamma$ the curve surrounding $S^q$. The two relations in \eqref{main} which define $S^q$ represent the implicit inequalities of two variables $x$ and $y$ that can be drawn with implicit functions plotting commands available in software such as Matlab (e.g.,  fimplicit), Mathematica (e.g.,  ContourPlot), or by using the free software Desmos \cite{desmos}.
In Figure  \ref{fig2}a, the case of $S^q$ for $q=0.5$ is presented.

\begin{remark}\label{rema}
As is known, if in the expression $x^y$,  $x$ is negative and $y$ is not an integer, the mathematical situation is somewhat ambiguous. With infinite numeric precision, the correct result of $x^y$ is mathematically well-defined without ambiguity. Certain values of $y$ yield an imaginary number as a result, while other values of $y$ result in a real-valued result. Specifically, when $x < 0$, the result of $x^y$ is real-valued exactly when $y$ can be written as a fraction, $m/n$, where $m$ is an integer and $n$ is an odd integer. Furthermore, the result is positive when $m$ is even but negative when $m$ is odd. When $y$ cannot be written as such, the result would be an imaginary number (see, e.g.,  \cite{wiki}). For example, in Matlab, for negative base $x$ and non-integer $y$, the power function returns complex results. A solution could be to use the function \emph{nthroot}. However, note that in IEEE floating-point computations
\[
x^y=\exp(y\log(x)).
\]
Moreover, the $\log$ function domain includes negative and complex numbers, which can lead to unexpected results if used unintentionally. Therefore, $x^y$ may be a complex number, if $x<0$ and $y$ is noninteger. In these cases, in some software, such as Matlab, imaginary parts of complex arguments are ignored to the detriment of the shape of $S^q$.
\end{remark}
To summarize, for those values $x$ and $y$ for which $E_a<-\pi/2$,  they are situated in quadrant III, $\cos(E)<0$ and, therefore, $\cos^q(E)$ is not well-defined . The consequence is that, in addition to the expected area, for some values of $q$, $S^q$ could present unexpected additional parts. Thus, for $q=0.8$, using Desmos, the domain $S^{0.8}$ is presented in Figure  \ref{fig2}b.

\begin{figure}
\centering
\includegraphics[scale=.7]{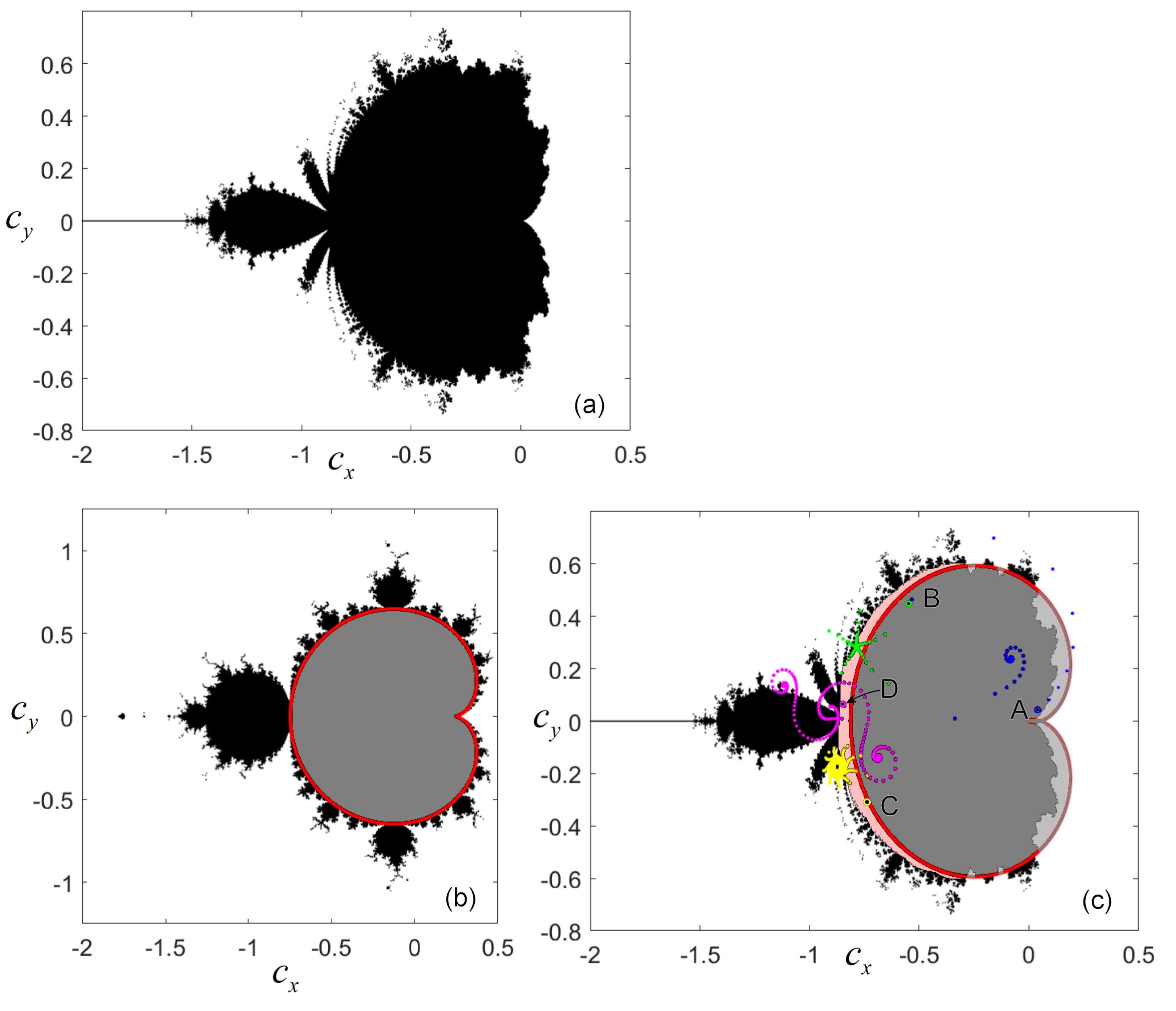}
\caption{Mandelbrot sets of IO and FO and their stability domains. (\textbf{a}) IOM set; (\textbf{b}) FOM set and the main cardioid (red plot) surrounding the stability domain for $q=0.85$.}
\label{fig3}
\end{figure}

In these cases, the supplementary domain can be determined via the condition $|\arg z|>q \pi/2$ (Matignon condition), for $q\in(0,1)$, and the lines $|\arg z|=q \pi/2$ being tangent to $S^q$. However, this is not possible in some more complicated cases (Figure  \ref{fig3}d, Section \ref{sec3}).


\vspace{-18pt}

\vspace{6mm}
To verify if $S^q$ in the considered cases of $q=0.5$ and $q=0.8$ contains the stability points, consider the points $X_i$, $i=1,...,4$ (images of the underlying complex numbers) for the case of $q=0.5$ and, $Y_i$, $i=1,2,...,7$ for the case of $q=0.8$, respectively (see Figure  \ref{fig2}), and denote
\begin{equation*}\label{aq}
A_q(z):=|z|-\left(2\cos\frac{|\arg z|-\pi}{2-q}\right)^q,~~~ q\in(0,1).
\end{equation*}
The graph of $A_q(z)$ is the boundary curve $\Gamma$ and
\[
S^q=\Bigg\{z\in \mathbb{C}: A_q(z)<0~~~\text{and}~~~ |\arg z|>\frac{q\pi}{2}\Bigg\}.
\]
Therefore, the images of $z$ for which $A_q(z)>0$, or $A_q(z)$ are complex values are not $S^q$ points.
The values of $A_q$ at the points $X_i$, and $Y_i$ are presented in Table \ref{tab1} and Table \ref{tab2}, respectively.

Regarding the numerical implementation, if in the atan2 function, $x$ and $y$ have different signs, then atan2 and arctan (atan) could have different values. Thus, if in some software the atan2 function is not implemented, it can be emulated with

\begin{equation*}\label{atan}
\text{atan2}(y,x)= \text{arctan}\frac{y}{x}+\frac{\pi}{2}\text{sign}(y)(1-\text{sign}(x)).
\end{equation*}

\begin{table}
\caption{Value
 of $A_{q}$ at points $X_i$, $i=1,...,4$, for $q=0.5$ (see also Figure  \ref{fig2}a). }\label{tab1}%
\setlength{\tabcolsep}{18.5pt}
\begin{tabular}{ccccc}
\toprule
$z$& $X_{1}\not\in S^{0.5}$ & $X_{2}\in S^{0.5}$ & $X_{3}\in S^{0.5}$ & $X_{4}\in cl(S^{0.5})$ \\
\midrule
$A_{0.5}(z)$ & 0.500 & $-$0.609 & $-$0.4142 & 0\\
 \bottomrule
\end{tabular}

\end{table}

\vspace{-10pt}

\begin{table}
\caption{Value of $A_{q}$ at points $Y_i$, $i=1,2,...,7$, for $q=0.8$ (see also Figure  \ref{fig2}b). }\label{tab2}%
\setlength{\tabcolsep}{4.7pt}
\begin{tabular}{cccccccc}
\toprule
$z$& $Y_{1}\not\in S^{0.8}$ & $Y_{2}\not\in S^{0.8}$ & $Y_{3}\not\in S^{0.8}$ & $Y_{4}\in S^{0.8}$ & $Y_{5}\in S^{0.8}$ & $Y_{6}\not\in S^{0.8}$ & $Y_{7}\in cl(S^{0.8})$ \\
\midrule$A_{0.8}(z)$ & $\in \mathbb{C}$ & $\in \mathbb{C}$ & $\in \mathbb{C}$ & $-0.5583$ & $-0.3823$ &
$0.5064$ & $1\times 10^{-4}$
\\%
 \bottomrule
\end{tabular}

\end{table}

\vspace{6mm}
An animation which shows the variation in the stability domain for $q\in(0,1)$, where several such cases appear, is presented as a supplementary video.

\section{Stability of the Mandelbrot Map of Fractional Order}\label{sec3}

Next, consider the fractional discretization of the Mandelbrot map

\begin{equation}\label{prima}
\Delta^q z(t)=f_c(z(t+q-1),~~z(0)=0,~~ q\in(0,1),~~t=0,1,2,...,
\end{equation}
where
\[
f_c(z(t+q-1)=z^2(t+q-1)+c,
\]
and $z,c\in \mathbb{C}$, with $c$ being the complex parameter. For the initial value problems for real fractional discrete systems, see, e.g.,  \cite{bibopt,bibnoua} ({in \cite{x4}, several properties of the complex Mandelbrot map of fractional order and the matlab code to draw the fractal are presented)}.

The numerical integral of the initial value problem \eqref{prima} is (see \cite{bibnoua} for the solution of real fractional-order systems)
\begin{equation}\label{ecus}
 z(n)=z(0)+ \frac{1}{\Gamma(q)}\sum_{i=1}^n\frac{\Gamma(n-i+q)}{\Gamma(n-i+1)}f_c(z(i-1)),~ z(0)=0, ~{n\in \mathbb{N}^*},
\end{equation}

For $q=0.85$, one obtains the fractional-order Mandelbrot (FOM) set, presented in Figure  \ref{fig3}b. The fractal is obtained with the FO\_mandelbrot.m code \cite{matlab}.

To obtain the stability domain of the stable fixed points (period 1) of the FOM map, consider next, for computational reasons, the problem in the Cartesian plane, with $c=c_x+ic_y$, $c_x,c_y\in \mathbb{R}$.
Contrary to the integer-order case, where the fixed points are found from the equation $f_c(z)=z$, here, the fixed points are obtained by solving the equation $f_c(z)=0$, i.e., $z^2+c=0$, with the solutions: $z^*_{1,2}=\pm i\sqrt c$. For the fixed point $z_1^*=i\sqrt{c}$, after some calculations, the eigenvalues are obtained as
\begin{equation}\label{eqq}
\lambda_{1,2}=\sqrt{2\sqrt{c_x^2+c_y^2}-2c_x}\pm i \sqrt{2\sqrt{c_x^2+c_y^2}+2c_x},
 \end{equation}
while, for $z_2^*=-i\sqrt{c}$,
\begin{equation}\label{eqq2}
\lambda_{3,4}=-\sqrt{2\sqrt{c_x^2+c_y^2}-2c_x}\pm i\sqrt{2\sqrt{c_x^2+c_y^2}+2c_x}.
 \end{equation}

Consider the eigenvalue $\lambda_1$, the reasoning for $\lambda_2$ being similar. Then, the useful parametric representation of $\Gamma$ in coordinates $c_x$ and $c_y$, defining $\Gamma(c)$, are determined from the following equations (see \eqref{para}):

\begin{equation*}\label{ecc}
\begin{array}{l}
\sqrt{2\sqrt{c_x^2+c_y^2}-2c_x}=-2^q\cos^q\theta \cos((2-q)\theta) \\
\sqrt{2\sqrt{c_x^2+c_y^2}+2c_x}=-2^q\cos^q\theta \sin((2-q)\theta),~~|\theta|\leq\pi/2.
\end{array}
\end{equation*}
The solutions, which represent the parametric equations of the curve $\Gamma$, are

\begin{equation}\label{crdiof}
\begin{array}{l}
c_x=-2^{2q-2}\cos ^{2q}\theta\cos(2\theta(q-2)),\\[.2cm]
c_y=2^{2q-1}\sin(\theta(2-q))\cos^{2q}\theta\cos(\theta(q-2)),~~|\theta|\leq\pi/2.
\end{array}
\end{equation}

In Figure  \ref{fig3}a, the integer-order Mandelbrot (IOM) set is drawn together with the main cardioid $\Gamma$ which has the equation (see, e.g.,  \cite{bibzece,bibunspe,bibdoispe})
\[
|1-\sqrt{1-4c}|=1.
\]
or, in parametric form,
\[
\begin{array}{l}
4c_x=2\cos\theta-\cos(2\theta),\\
4c_y=2\sin\theta-\sin(2\theta),~~|\theta|\leq \pi.
\end{array}
\]

The interior of the curve $\Gamma$ corresponds to the stable fixed points, obtained from the equation $f_c(z)=z$, i.e. $z_{1,2}^*=(1\pm \sqrt{1-4c})/2$.

The curve $\Gamma$ for the FOM set for $q=0.85$, which surrounds the stability domain of the fixed point $z_1^*$, $S^{0.85}$, is presented in Figure \ref{fig3}c. Note that in both the IO case (Figure  \ref{fig3}a) and the FO case (Figure  \ref{fig3}b), the Mandelbrot sets (defined on the parametric plane of $c$ in coordinates $(c_x,c_y)$) and the stability domains (defined on the eigenvalues $\lambda$ in coordinates $(\lambda_x,\lambda_y)$) are overplotted.

As detailed in Section \ref{sec2}, to a fixed point $z^*$, the eigenvalues $\lambda$ given by \eqref{eqq} or \eqref{eqq2} correspond. For a point $\lambda\in \Gamma$, $z*$ loses its stability, while if $\lambda\notin S^{0.85}$,  $z^*$ is unstable. If $z^*$ is asymptotically stable, $\lambda$ belongs to $S^{0.85}$, and reversely, a point $\lambda$ within $S^{0.85}$ corresponds to a point $c$ for which $z^*$ is asymptotically stable.

Let us  numerically verify this property. Consider $\lambda=-0.5701 +i0.3019\in S^{0.85}$ (magenta point $1$ of coordinates $(-0.5701,0.3019)$ in Figure  \ref{fig3}b). Because $\Re{\lambda}<0$ and $\Im{\lambda}>0$, one chooses \eqref{eqq2} with right-hand side having signs $-$ and $+$, respectively. To find the corresponding point $c$ of this value of $\lambda$, one has to solve the problem \eqref{eqq} with unknown $c_x$ and $c_y$. The solutions are $c_x=-0.0585$ and $c_y=0.0861$ (magenta point 2 in Figure  \ref{fig3}b. To verify the stability of the fixed point corresponding to this value of $c$, one integrates the system with \eqref{ecus}. The orbit starting from point $c$ tends to be point 3 of the coordinates $(-0.2845, 0.1510)$ (Figure  \ref{fig3}b), which is an approximation with an error of $1\times 10^{-3}$ of the fixed stable point $i\sqrt{c}=i\sqrt{-0.0585+0.0861i}=-0.2851 + 0.1510i$.

Interestingly, in comparison with the IO case (Figure  \ref{fig3}a), $S^{0.85}$ does not fill the entire FOM set (gray fill in \mbox{Figure  \ref{fig3}b).} However, points within this gray area of $S^{0.85}$, and similarly for all studied cases of $S^q$, lead to the same result: the entire domain $S^q$ contains eigenvalues points for which the fixed points are asymptotically stable. Consider, for example, the eigenvalue $\lambda=0.1464 -i0.2268$ (green point 1 of coordinates $(0.1464, -0.2268)$ in \mbox{Figure  \ref{fig3}b).} By solving system \eqref{eqq}, one obtains the point $c_x=0.0075$ and $c_y=0.0166$ (green point 2 in Figure  \ref{fig3}b), from which the integral \eqref{ecus} generates the orbit tending to point 3 of the coordinates $(-0.0725, 0.1134)$. This, with a precision $1\times10^{-3}$, represents the approximation of the fixed point $i\sqrt{c}=-0.0732 + 0.1134i$.

On the other hand, points outside of $S^{0.85}$ are eigenvalues for which the underlying points $c$ do not belong to the FOM set and, therefore, their fixed points are not stable. Consider, for example, point 1 (blue plot in Figure  \ref{fig3}b) of coordinates $\lambda_x=0.1231$ and $\lambda_y= 0.5590$. Solving system \eqref{eqq}, one obtains $c_x=0.0743$, and $c_y=0.0344$ (point 2 blue outside the FOM set in Figure  \ref{fig3}b) from where the orbit diverges in agreement with the definition of the Mandelbrot set of IO or FO.

Among several differences and resemblances between the FOM set and the IOM \mbox{set \cite{sub},} one can see that the origin of the so-called ``Elephant Valley'' (see, e.g.,  \cite{bibzece}) in the case of the IOM set is located at $(1/4,0)$, while for the FOM set, it is located at $(0,0)$  (see dotted line in Figure \ref{fig3}).

\begin{figure}
\includegraphics[scale=.6]{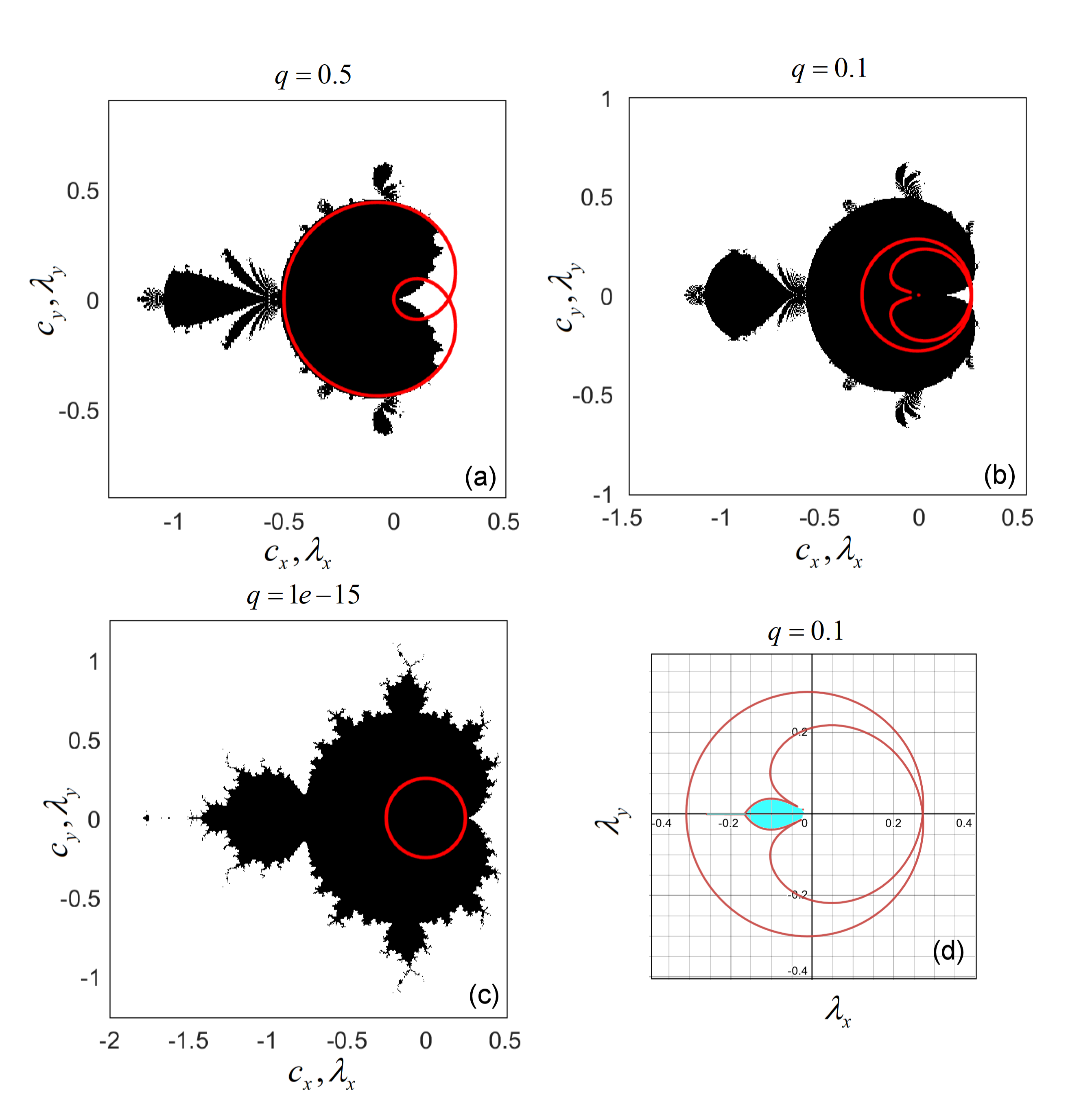}
\caption{Mandelbrot sets of FO and their stability domains for three extra cases. (\textbf{a}) The FOM set and $S^q$ for the limit case  $q=0.5$; (\textbf{b}) FOM set and $S^q$ for $q=0.1$; (\textbf{c}) FOM set and $S^q$ for $q=1\times 10^{-15}$; (\textbf{d}) $S^q$ for $q=0.1$, obtained in Desmos software.}
\label{fig4}
\end{figure}

While for $q\geq0.5$, $S^q$ fills the main body of the FOM set well enough, for $q<0.5$, the shape of $S^q$ does not cover the main body. In Figure \ref{fig4}a--c, for $q=0.5$, $q=0.1$ and $q=1\times10^{-15}$,
 the numerically near $0$ value,  $q=1\times10^{-15}$, is chosen instead of $\Gamma(0)$, which is not defined. In this case, in order to obtain more fractal details (which look similar to the IOM set), only $30$ iterations were used, compared to $1000$ iterations for the other cases). One can see that $q=0.5$ is the ultimate value for $q$. Thus, for $q$ values below this limit, $S^q$ is shrinking with respect to the main body of the FOM set. Another characteristic is that the shape of $\Gamma$ obtained with Desmos (Figure  \ref{fig4}d) presents an additional part (see also Figure  \ref{fig1}b). This image was obtained using both representations: the representation \eqref{main} and the parametric form \eqref{crdiof}.
\vspace{-14pt}

\vspace{6mm}
Therefore, the following conjecture is numerically sustained (see also \cite{fo}):

\begin{conjecture}\label{conju}
 For $q<0.5$, the surface of the stability domain $S^q$ of the fixed points of the FOM shrinks compared to the main body of the fractal set.
\end{conjecture}

\begin{remark}
After the first result on the non-periodicity of non-constant solutions of FO, continuous-time systems appeared \cite{tava}; this intriguing result has been extended to FO discrete systems too (see, e.g.,  \cite{mich}). Thus, contrary to the IOM set, the head contains points generating two-period stable cycles, and where the bulbs contain the points, generating multiple stable cycles of period 3.4, and so on (Figure  \ref{fig3}a). In the case of the FOM set, these apparent periodic orbits do not exist! For example, points in the head (yellow fill in Figure  \ref{fig3}b) do not generate periodic orbits. The same situations happen to all bulbs.
\end{remark}

\section{Conclusions}
In this paper, it is shown that drawing the stability domain $S^q$ for empiric values of $q\in(0,1)$ encounters the problem of the power of negative numbers which could be differently treated by the software used. Therefore, depending on the software used, the graphical results of the stability domain of the linear difference equation \eqref{unu} could be unexpected. The power function appearing in $S^q$ is applied to negative values when either the logarithm is a complex number or one obtains an unexpected result (see Remark \ref{rema}). In addition, the stability domain for the fixed points of an FOM map is determined. The asymptotical stability of the fixed points of the FOM map is verified numerically for $q>0.5$. For $q<0.5$, the shape of $S^q$ seems not to verify the filling characteristic encountered on the IOM set,  {a fact that leads to Conjecture \ref{conju}}.
Besides several differences between the Mandlebrot sets of IO and FO, such as the shape, the position, and the shape of the stability set, in the FO case, the head of the FOM set does not represent the points that generate stable periodic orbits.

\vspace{6pt}

References

\end{document}